\documentclass[10pt]{article}
\usepackage{graphicx}
\usepackage{amssymb}
\usepackage{epstopdf}
\usepackage{amsmath}
\usepackage{amsfonts}
\usepackage{color}
\usepackage[toc,page, title]{appendix}

\def\phi{{\varphi}}

\DeclareSymbolFont{AMSb}{U}{msb}{m}{n}
\DeclareMathSymbol{\N}{\mathbin}{AMSb}{"4E}
\DeclareMathSymbol{\Z}{\mathbin}{AMSb}{"5A}
\DeclareMathSymbol{\R}{\mathbin}{AMSb}{"52}
\DeclareMathSymbol{\Q}{\mathbin}{AMSb}{"51}
\DeclareMathSymbol{\I}{\mathbin}{AMSb}{"49}
\DeclareMathSymbol{\C}{\mathbin}{AMSb}{"43}

\def\be{\begin{equation}}
\def\ee{\end{equation}}
\def\ber{\begin{eqnarray}}
\def\eer{\end{eqnarray}}

\def\beq{\begin{equation}}
\def\eeq{\end{equation}}

\begin{document}

\addtolength{\textheight}{0 cm} \addtolength{\hoffset}{0 cm}
\addtolength{\textwidth}{0 cm} \addtolength{\voffset}{0 cm}

\newenvironment{acknowledgement}{\noindent\textbf{Acknowledgement.}\em}{}

\setcounter{secnumdepth}{5}

\newtheorem{theorem}{Theorem}[section]
\newtheorem{lemma}[theorem]{Lemma}
\newtheorem{coro}[theorem]{Corollary}
\newtheorem{remark}[theorem]{Remark}
\newtheorem{claim}[theorem]{Claim}
\newtheorem{conj}[theorem]{Conjecture}
\newtheorem{definition}[theorem]{Definition}
\newtheorem{application}{Application}
\newtheorem{example}[theorem]{Example}
\newtheorem{corollary}[theorem]{Corollary}
 \newtheorem{proposition}[theorem]{Proposition}

\title{Solutions to multi-marginal optimal transport problems concentrated  on several graphs\footnote{  B.P. is pleased to acknowledge the support of a University of Alberta start-up grant and National Sciences and Engineering Research Council of Canada Discovery Grant number 412779-2012.  A. M. greatly acknowledges  the support of  a grant from the Natural Sciences and Engineering Research Council of Canada (417885-2013).
}}  
\author{Abbas Moameni \footnote{School of Mathematics and Statistics,
Carleton University,
Ottawa, Ontario, Canada K1S 5B6, 
momeni@math.carleton.ca} and Brendan Pass\footnote{Department of Mathematical and Statistical Sciences, 632 CAB, University of Alberta, Edmonton, Alberta, Canada, T6G 2G1 pass@ualberta.ca.}}

\maketitle

\begin{abstract}
We study solutions to the multi-marginal Monge-Kantorovich problem which are concentrated on several graphs over the first marginal.  We first present two general conditions on the cost function which ensure, respectively, that any solution must concentrate on either finitely many or countably many graphs.   We show that local differential conditions on the cost, known to imply local $d$-rectifiability of the solution, are sufficient to imply  a local version of the first of our conditions.  We exhibit two examples of cost functions satisfying our conditions, including the Coulomb cost from density functional theory in one dimension.  We also prove a number of results relating to the uniqueness and extremality of optimal measures.  These include a sufficient condition on a collection of graphs for any competitor in the Monge-Kantorovich problem concentrated on them to be extremal, and a general negative result, which shows that when the problem is symmetric with respect to permutations of the variables, uniqueness cannot occur except under very special circumstances.
\end{abstract}
.
\section{Introduction}

In this paper, we consider the multi-marginal Monge-Kantorovich transport problem.
Given Borel probability measures $\mu_1,..., \mu_n$ (the marginals) on smooth manifolds    $X_1,..., X_n,$ respectively, let  $\Pi(\mu_1,...,\mu_n)$ be  the set of Borel probability measures on $X_1 \times X_2 \times... \times X_n$  which have $X_i$-marginal $\mu_i$ for each $i \in \{1,2,...,n\}.$  Letting \[c: X_1 \times X_2 \times... \times X_n  \to \R,\]  be  a measurable  cost function, the transport cost associated to a transport plan
 $\pi \in \Pi(\mu_1,...,\mu_n)$ is defined by
\[I_c(\pi)=\int_{X_1 \times X_2 \times... \times X_n} c(x_1,...,x_n) \, d \pi.\]
The multi-marginal Monge-Kantorovich transport problem is then the following minimization problem:
\[  \inf \Big\{ I_c(\pi); \pi \in \Pi(\mu_1,...,\mu_n)\Big\}. \qquad \qquad (MK)        \]

Under quite mild conditions on $c$ and the $\mu_i$, it is well known that $ (MK)$ admits a minimizer.  Our purpose here is to study the geometric structure of optimal measures $\pi \in  \Pi(\mu_1,...,\mu_n)$.

Note that when $n=2$,  $(MK)$ is precisely the classical optimal transport problem of Monge and Kantorovich.  In this case, assuming that the first marginal $\mu_1$ is absolutely continuous with respect to local coordinates, the celebrated \emph{twist} condition implies that the solution is concentrated on a graph over $x_1$, and is unique \cite{L,G95,GM96,Caf,V}.  While this condition holds for many important cost functions, there are also a variety where it fails, including \emph{any} smooth cost on a compact manifold.  Recently, there have been some developments on two marginal problems where the twist condition fails \cite{B-C, C-M-N, Moa}.  In particular, one of the present authors developed two conditions, known as the \emph{$m-$twist} and \emph{generalized twist} conditions, which ensure respectively that the solution concentrates on several or countably  many (rather than one) graphs \cite{Moa}.  A related line of recent research concerns sufficient conditions for uniqueness of the optimal measure, one interesting goal being to exhibit smooth costs on arbitrary manifolds for which optimal plans are unique, despite the fact that such plans are not generally concentrated on graphs. The interested reader is referred to a recent manuscript \cite{M-R} and also to \cite{A-K-M, B-C, C-M-N, G-M2} for more examples.  We recall that a measure $\gamma$ in the convex and weakly compact set $\Pi(\mu_1, \mu_2)$ is called \textit{extremal} if $\gamma=t \gamma_1 +(1-t)\gamma_2$ with $\gamma_1, \gamma_2 \in \Pi(\mu_1, \mu_2)$ and $0<t<1$ imply that $\gamma_1=\gamma_2.$  It is well known that if $(MK)$ admits a unique solution, then this solution must be extremal in $\Pi(\mu_1, \mu_2)$; an important aspect of the research on uniqueness is the characterization of extremal points of $\Pi(\mu_1, \mu_2)$.  It is  known that the extremal points of $\Pi(\mu_1, \mu_2)$ are not necessarily  concentrated on a single graph. Indeed, any extremal point of $\Pi(\mu_1, \mu_2)$ is known to concentrate on the union of the   graph of a function from $X_1$ to $X_2$ and the  antigraph of a function from $X_2$ to $X_1.$  We refer the interested reader to \cite{Moa3} where a sufficient  and almost necessary  condition for the structure of the support of extremal measures in $\Pi(\mu_1, \mu_2)$ is established.

Interest in the $n\geq 3$ case (also known as the multi-marginal case) has picked up recently.   Here, an  analogue of the twist condition, known as the \emph{twist on splitting sets} condition  has recently been developed; this condition ensures that the solution to the multi-marginal problem concentrates on a  graph over the first variable $x_1$ (and is unique) \cite{K-P}.   In contrast to the two marginal setting, where the twist condition is reasonably mild (at least for spaces with trivial topology), the twist on splitting sets condition is very strong and there are many important examples for which it fails.  In particular, examples arising in density functional theory in physics \cite{B-D-G, C-F-K} and  roommate type problems in economics \cite{cgs} involve equal marginals and cost functions which are symmetric under all permutations of the arguments  (see also \cite{G-M2}).  For problems with this symmetric structure, it is known that there are always solutions that do not concentrate on a single graph (except in the very simple special case when the diagonal $\{(x,x,...x)\}$ supports a unique optimizer) \cite{P3}.  In this paper, we focus on what can be  said in the multi-marginal case when the solution does not concentrate on a single graph; in particular, we investigate the possibility that the solution may be concentrated on several graphs. 
 
Specifically, we consider the $m-$twist condition  on splitting sets, and the related generalized twist on splitting sets condition, (see Definition \ref{mtwist}) which are analogues of the conditions in \cite{Moa} in the two marginal case.    These conditions were recently introduced by one of the authors  in  \cite{Moa2}, and assuming the first marginal is absolutely continuous with respect to local coordinates, they imply, respectively, that any optimal plan $\pi$ concentrates on at most $m$ graphs, or countably many graphs, over the first variable.

The  $m-$twist condition  on splitting sets is somewhat more flexible than the twist on splitting sets condition and, in particular, it seems relevant to symmetric costs which arise in the applications outlined above.  We show, for instance, that the Coulomb cost (see \eqref{coulombcost} below) from density functional theory is $m$-twisted on splitting sets in one dimension.  Our work therefore offers a different perspective on a recent result of Colombo, De Pascale and Di Marino \cite{cdpdm}, who constructed explicit solutions to $(MK)$ for this cost, when the marginals are all equal; their solutions concentrated on $(n-1)!$ graphs (where $n$ is the number of marginals).  In fact, our result also applies to the case when the marginals are not equal; on the other hand, we are able to prove only that solutions concentrate on at most $n!-(n-1)!$ graphs.  It is unclear to us whether there are in fact marginals for which our upper bound on the number of graphs is sharp; if so, the marginals cannot all be equal.

We also demonstrate that a local differential condition on $c$, which is known to imply that the optimizer concentrates on a low dimensional subset of the product space, also implies a local version of the twist condition  on splitting sets. This local twist on splitting sets in turn implies the generalized twist on splitting sets, so  that each optimal plan concentrates on the union of the graphs of  at most a countable number of maps.

Following the recent developments on uniqueness issues for non graphical solutions in the two marginal case, we  establish here some results which can be useful to derive uniqueness of  solutions concentrated on several graphs in the multi-marginal setting, and exhibit some examples.  We also prove a general, negative result, implying that the permutation symmetric problems described above can have unique solutions only in very special circumstances.  For solutions concentrated on several graphs, however, these circumstances are somewhat more flexible than in the single graph case (where only the trivial, diagonal solution $(Id,Id,...Id)_\# \mu$ has the possibility of being a unique optimizer); we provide an example of a symmetric three marginal problem where a unique solution concentrates on two graphs and at least one of the components of each is not the identity mapping.

The next section is devoted to preliminary definitions.  In the third section, we state our results for costs which satisfy the $m$ and generalized twist on splitting sets conditions.  We then prove that the local differential conditions alluded to above imply  local twist on splitting sets, and outline several examples.  In the  fourth section, we prove several results related to uniqueness of solutions concentrated on multiple graphs.  We show that under the generalized twist on splitting sets condition, a collection of several graphs containing the support of every optimizer is uniquely determined.  We then establish conditions under which any solution concentrated on several graphs must be extremal, and prove that solutions to symmetric problems cannot be unique when $n \geq 3$, except under very particular conditions.

\section{Definitions and preliminaries}
Here we recall several definitions which will be needed throughout the paper.  We begin with the definition of $c$-splitting sets, modified slightly from \cite{K-P}.
\begin{definition} A set $S \subset  X_1 \times X_2\times ...\times X_n$ is a $c$-splitting set (or simply a splitting set) if there exist
Borel functions $u_i : X_i \to [-\infty,\infty)$  such that for all  $(x_1, x_2, ..., x_n),$
\begin{equation}\label{eqn: splitting inequality}
 \sum_{i=1}^n  u_i(x_i) \leq  c(x_1, x_2, ..., x_n)
\end{equation}
with equality whenever $(x_1, x_2, ..., x_n) \in S$.  The $n$-tuple $(u_1,...,u_n)$ is called the $c$-splitting tuple for $S.$
\end{definition}
 
The relevance of splitting sets to $(MK)$ is that, as can be easily seen from the dual formulation found in \cite{Ke}, any optimal measure is concentrated on a $c$-splitting set.
 
Next, we recall the following relaxations, introduced in \cite{Moa2}, of the twist on splitting sets conditions from \cite{K-P}.  In what follows, $D_{x_1}c(x_1,x_2,...,x_n)$ denotes the differential of $c$ with respect to the first variable $x_1$.
\begin{definition}\label{mtwist} Let  $c$ be a Borel measurable function and denote by $D_1(c)$ the set of points at which $c$ is  differentiable with respect to the first variable. \\
1.  {{\bf $\mathbf{m}$-twist condition:}} 
We say that 
$c$ is $m$-twisted on splitting sets 
if  for any $c$-splitting set  $S \subset  X_1 \times X_2 \times  ...\times X_n$ , any $x_1 \in X_1$ and  any $p$ in the cotangent space $T^*_{x_1}X_1$ to $X_1$ at $x_1$,
 the cardinality of the set 
\[\Big\{( x_1,x_2, ..., x_n) \in S\cap D_1(c);\, D_{x_1}c( x_1,  x_2, ...,  x_n)=p \Big\},\]
is at most $m.$  We  say that $c$  locally satisfies the $m$-twist condition  on splitting sets  if for any  $(x_1,  x_2, ...,  x_n)\in S\cap D_1(c)$ 
there exists a neighborhood $U \subseteq X_2 \times...\times X_n$ of $(x_2,...,x_n)$ such that
 the cardinality of the set \[\Big\{y \in U,\, (x_1,y) \in S \cap D_1(c), \, D_{x_1} c(x_1,y)=D_{x_1} c(x_1,x_2,...,x_n)\Big\},\]
is at most $m.$\\
2. {{\bf Generalized-twist condition:}} We say that 
$c$ satisfies the generalized twist condition on splitting sets if for any
 $c$-splitting set  $S \subset  X_1 \times X_2 \times  ...\times X_n$, any $x_1 \in X_1$ and  any $p$ in the cotangent space $T^*_{x_1}X_1$ to $X_1$ at $x_1$, the set \[\Big\{( x_1,x_2, ..., x_n) \in S\cap D_1(c);\, D_{x_1}c(x_1,  x_2, ...,  x_n)=p \Big\},\]
is a finite subset of $S$.
\end{definition}
\begin{remark}\label{rem: cont splitting functions}
Given a $c$-splitting set $S$ and splitting functions $u_1,....u_n$, it is well known that one can construct new splitting functions, $\bar u_1,....\bar u_n$ for $S$ such that, for each $i=1,2....,n$, we have
\begin{equation}\label{eqn: c-conjugate}
\bar u_i(x_i) =\inf_{x_j, j\neq i}c(x_1,...x_n) -\sum_{j\neq i }\bar u_j(x_j).
\end{equation}
Assuming $c$ is globally Lipschitz, it is well known that each $u_i$ is as well \cite{McCann01} (unless $u_i$ is identically $-\infty$, in which case the splitting set is empty).  Therefore, the set $ \{ c(x_1,...x_n) =\sum_{ i=1 }^n\bar u_i(x_i)\}$, which contains $S$, is closed (and hence compact if each $X_i$ is).  
\end{remark}

The following result provides a connection between the generalized twist condition on splitting sets and local $m$-twistedness on splitting sets.

\begin{proposition}\label{non}
Assume that $c$ is continuously differentiable with respect to the first variable and each $X_i$ is compact. 
 If, for some $m \in \mathbb{N},$ the function  $c$ is locally $m$-twisted on splitting sets, then $c$ satisfies the generalized-twist condition on splitting sets.  
\end{proposition}

\textbf{Proof}. Assume that $S \subset X_1 \times ...\times X_n$ is a $c$-splitting set.  

Fix $(\bar x_1, ...,\bar x_n) \in S$. We need to show that
the set \[L=\Big\{(\bar x_1, x_2,...,x_n) \in S; \, D_1 c(\bar x_1 , \bar x_2,...,\bar x_n)=D_1 c(\bar x_1, x_2,...x_n)\Big\},\]
is finite.   By  Remark \ref{rem: cont splitting functions}, we can find a compact splitting set, $\tilde S$, such that $S \subseteq \tilde S$; we will in fact prove the stronger statement, that
 \[\tilde L=\Big\{(\bar x_1, x_2,...,x_n) \in \tilde S; \, D_1 c(\bar x_1 , \bar x_2,...,\bar x_n)=D_1 c(\bar x_1, x_2,...x_n)\Big\},\]
is finite (note that $L \subseteq \tilde L$).
If $\tilde L$ is not finite there exists a countably infinite subset $\{(\bar x_1, x^k_2,...x^k_n)\}_{k \in \mathbb{N}} \subset\tilde L.$  
Since $\tilde S$ is  compact  then the sequence  $\{(\bar x_1, x^k_2,...x^k_n)\}_{k \in \mathbb{N}}$ has an
accumulation point $(\bar x_1, x^0_2,...x^0_n)\in S$ and there exists a subsequence still denoted by $\{(\bar x_1, x^k_2,...x^k_n)\}_{k \in \mathbb{N}}$ 
such that $x_i^k \to x_i^0$ as $k \to \infty$ for $i=2,...,n.$  Since
 $D_1 c$ is continuous it follows that
$(\bar x_1, x^0_2,...x^0_n) \in \tilde L.$ Since $c$ is  locally $m$-twisted on $\tilde S$, this leads to a contradiction   as $(\bar x_1, x^0_2,...x^0_n) $ is an accumulation point of the
sequence $\{(\bar x_1, x^k_2,...x^k_n)\}_{k \in \mathbb{N}}$ and 
\[D_1 c(\bar x_1, x^0_2,...x^0_n)=D_1 c(\bar x_1, x^k_2,...x^k_n), \qquad \forall  k \in \mathbb{N}.\] This completes the proof. \hfill $\square$

We recall next an important property of splitting sets, which will be useful at various points throughout the paper.  

\begin{definition}\label{def: cmonotonicity}
A set $S \subset X_1 \times X_2 \times...\times X_n$ is \emph{$c$-monotone} if the following holds.  Whenever $x, \bar x \in S$, and $P_+,P_-$ are two non empty disjoint sets of indices such that $P_+ \cup P_-=\{1,2,...,n \}$, we have
$$
c(x) +c(\bar x) \leq c(x_+, \bar x_-)+c(\bar x_+,  x_-).
$$
Here, we have decomposed $x =(x_+,x_-)$ and $\bar x =(\bar x_+,\bar x_-)$ in the obvious way; that is, $x_+ =(x_i)_{i\in P_+}$, with analagous definitions for $x_- ,\bar x_+$ and $\bar x_-$. 
\end{definition}
It is well known that any $c$-splitting set is $c$-monotone (see \cite{P2}\cite{K-P})).

Finally, we define precisely some notation describing measures concentrated on several graphs.
\begin{definition}\label{union} Let $X$ and $Y$ be Polish spaces with Borel probability measures $\mu$ on $X$ and $\nu$ on $Y.$
 We say that a measure $\gamma \in \Pi(\mu, \nu)$ is concentrated on the graphs of  measurable maps  $\{T_i\}_{i=1}^k$  from $X$ to $Y$,
if there exists a sequence of measurable non-negative  functions
$\{\alpha_i\}_{i=1}^k$ from $X$ to $\R$ with $\sum_{i=1}^k \alpha_i(x)=1$ ($\mu$-almost surely) such that for each measurable set $S \subset X \times Y,$
\[\gamma(S)=\sum_{i=1}^k\int_X \alpha_i(x) \chi_S(x,T_ix) \, d\mu,\]
where $\chi_S$ is the indicator function of the set $S.$ In this case we write $\gamma=\sum_{i=1}^k \alpha_i(Id \times T_i)_\# \mu.$
\end{definition}

\section{Multi-graph solutions}

The following characterization and its proof can be found in \cite{Moa2}; it asserts that, under the $m-$twist on splitting sets condition, any solution to the Kantorovich problem concentrates on (at most) $m$ graphs over the first variable.

\begin{theorem} \label{main} Assume that  the cost function $c$   satisfies  the $m-$twist condition on splitting sets, 
$\mu_1$ is non-atomic and  any function $u_1(x_1)$ that is of form \eqref{eqn: c-conjugate} and not identically infinite is differentiable $\mu_1$-almost surely on its domain.
  Then for each optimal plan
  $\gamma $ of  $(MK)$ with $Supp(\gamma) \subset D_1(c),$
 there exist $k \leq m,$ a sequence of non-negative measurable real functions 
$\{\alpha_i\}_{i=1}^k $ on $X_1$ 
 and,  Borel measurable maps $G_1,...,G_k: X_1 \to  X_2 \times ...\times X_n$ such that
\begin{equation}\label{late}\gamma =\sum_{i=1}^k \alpha_i (\text{Id} \times  G_i)_\# \mu, \qquad \end{equation}
where  $\sum_{i=1}^k \alpha_i(x) =1$ for $\mu_1$-a.e. $x \in X_1$.\\
Moreover, if one replaces the $m$-twist condition by the generalized-twist condition then   each optimal plan
  $\gamma $ of  $(MK)$ is of the form (\ref{late}) for some  $k \in \mathbb{N}\cup \{\infty \}.$
\end{theorem}
Although this result asserts quite strong conclusions about the structure of optimizers, the $m-$ twist on splitting sets (respectively, generalized twist on splitting sets) hypothesis seems quite difficult to verify.  Below, we establish that certain local differential conditions on the cost actually suffice for the generalized twist on splitting sets condition.

\subsection{A differential condition for local twist on splitting sets}

Here we establish a differential condition on $c$ which is sufficient for the local $1-$twist on splitting sets condition.  The condition first appeared in \cite{P2}, where is was shown to imply that the solutions concentrate  on low dimensional subsets of the product space.
Assume that each $X_i$ is a $d$-dimensional smooth manifold.  As in \cite{P2}, we let $g$ be the off diagonal part of the Hessian of $c$; that is, in block form,

\begin{equation}\label{gmatrix} \qquad
g=
\begin{bmatrix}
 0 &D^2_{x_1x_2}c &...&...&D^2_{x_1x_n}c \\
 D^2_{x_2x_1}c &0 & D^2_{x_2x_3}c&...&D^2_{x_2x_n}c \\
...& ...&...&...&...\\
 D^2_{x_nx_1}c & D^2_{x_nx_2}c&... &...&0\\
\end{bmatrix},
\end{equation}
where each $D^2_{x_ix_j}c =(\frac{\partial ^2c}{\partial x_i^\alpha x_j^\beta})_{\alpha\beta}$, for $i \neq j$, is the $d \times d$ matrix of mixed, second order partial derivatives
We recall that the signature of a symmetric matrix is the ordered triplet, $(\lambda_+,\lambda_-,\lambda_0)$, where $\lambda_+$, $\lambda_-$ and $\lambda_0$ denote respectively the number of positive and negative eigenvalues, and the multiplicity of the zero eigenvalue.

We then have the following.
\begin{proposition}\label{gsignature} Let $c$ be twice differentiable. The following assertions hold:
\begin{enumerate}\item If  the symmetric, $mn \times mn$ matrix $g$ has signature $(nd-d,d,0)$ then $c$ is locally $1-$twisted on splitting sets.
\item If in addition each $X_i$ is compact then $c$ satisfies the generalized twist condition on splitting sets.
\end{enumerate}
\end{proposition}

\textbf{Proof.} We first prove assertion 1.  The proof is by contradiction; to this end, assume that $c$ is not locally $1-$twisted on splitting sets.  It follows that there exists  a splitting set $S$ and a point $(\bar x_1, \bar x_2, ..., \bar x_n) \in S$ such  that the set 

$$
\Big\{(x_2, ..., x_n):(\bar x_1,x_2, ..., x_n) \in S,\,D_{x_1}c(\bar x_1, \bar x_2, ..., \bar x_n)= D_{x_1}c(\bar x_1,  x_2, ...,  x_n) \Big\}
$$ 
intersects any open neighbourhood of $( \bar x_2, ..., \bar x_n)$ (at a point other than $( \bar x_2, ..., \bar x_n)$).
We can therefore take a sequence  $(x^k_2, ..., x^k_n)$ in this set converging to $(\bar x_2, ...,\bar  x_n)$, such that $(x^k_2, ..., x^k_n) \neq (\bar x_2, ..., \bar x_n)$ for all $k$.

Set 

\[v^k =\frac{(x^k_2, ..., x^k_n)-(\bar x_2, ..., \bar x_n)}{|(x^k_2, ..., x^k_n)-(\bar x_2, ..., \bar x_n)|}.\]  
As each $v^k$ has unit norm, we may pass to a convergent subsequence, so that $v^k \rightarrow v=(v_2,...,v_n) \neq 0$.  Note that the vector $(0,v) =(0,v_2,...,v_n)$ is then tangent to $S$ at $(\bar x_1, \bar x_2, ..., \bar x_n)$.

Furthermore,  the splitting set $S$ is necessarily $c-$monotone, 

and so, by a result in \cite{P2} we have

$$
(0,v)^T\cdot g \cdot (0,v) \leq 0
$$
This is equivalent to

$$
\sum_{i,j =2, i\neq j}^n v_i \cdot D^2_{x_ix_j}c \cdot v_j \leq 0.
$$
Now, since for each $k$ we have

$$
\frac{D_{x_1}c(\bar x_1, \bar x_2, ..., \bar x_n)- D_{x_1}c(\bar x^k_1,  x^k_2, ...,  x^k_n)}{|(x^k_2, ..., x^k_n)-(\bar x_2, ..., \bar x_n)|} =0,
$$
taking the limit as $k \rightarrow \infty$ yields
$$
0=\sum_{j=2}^nD^2_{x_1x_j}c\cdot v_j.
$$

But now let $w \in T_{\bar x_1}X_1$, $s \in \mathbb{R}$, and consider the vector $u = (w,0,0,...,0)+s(0,v)$.  We then have

$$
u^t\cdot g \cdot u =  \sum_{i,j=1, i\neq j}^k    u_i \cdot D^2_{x_ix_j}c \cdot u_j     =2sw\cdot \sum_{i=2}^nD^2_{x_1x_j}c\cdot \cdot v_j +\sum_{i,j =2, i\neq j}^n v_i \cdot D^2_{x_ix_j}c \cdot v_j =   s^2 v^t\cdot g \cdot v \leq 0.
$$

We have therefore found a $d+1$ dimensional vector space, $T_{\bar x_1}X_1 \oplus \text{span}(v) \subseteq T_{\bar x_1 \times \bar x_2,...,\bar x_n}X_1 \times X_2 \times ...X_n$,  on which $g$ is negative definite, contradicting the assumption that $g$ has only $d$ timelike directions, and completing the proof of assertion 1.
Proposition \ref{non} then immediately implies the second assertion. $\square$\\

The last proposition together with Theorem \ref{main} imply that each optimal plan is concentrated on the union of the graphs of at most a countable number of maps.

\subsection{Application: the one dimensional  Coulomb cost with $n$ marginals.}
Let $c: \R^n \to \R$ be the one dimensional Coulomb cost; that is, take each $X_i =\mathbb{R}$ and set

\begin{equation}\label{coulombcost}
c(x_1,...,x_n)=\sum_{1\leq i<j\leq n} \frac{1}{|x_i-x_j|}, \qquad \forall (x_1,...,x_n)\in \R^n.
\end{equation}

This cost has important applications in density functional theory in physics \cite{C-F-K}\cite{B-D-G}.  A recent paper by Colombo, De Pascale and  Di Marino \cite{cdpdm} solves the optimal transport problem with this cost explicitly when the marginals $\mu_i$ are all the same, which is the most important case in view of physical applications.  The solutions they exhibit all concentrate on several graphs.  Here, we show that this cost satisfies the $m$ twist on splitting sets condition, which, by virtue of Theorem \ref{main}, uncovers a new perspective on their result, as well as showing that the multi-graph solution structure persists to the setting where the marginals  differ.
\begin{proposition}\label{col} The one dimensional Colulomb cost \eqref{coulombcost} is $n! -(n-1)!$ twisted on splitting sets.
\end{proposition}
\textbf{Proof.}
Let $S$ be a splitting set and $u_1,...,u_n$ the corresponding splitting functions.  Fix $x_1 \in \mathbb{R}$ and $p\in \mathbb{R}$.  As the splitting functions never take on the value $+\infty$, the inequality \eqref{eqn: splitting inequality} implies that there is no point $x_1,...x_n \in S$ with $x_i =x_j$ for some $i \neq j$.  We first prove that, for each permutation $\sigma$ on $n$ letters, there is at most one $(x_2,...,x_n) \in X_2\times...\times X_n$, such that $D_{x_1}c(x_1,...x_n) =p$ and 
$$
 (x_1,x_2,...x_m) \in S \cap A_{\sigma}
$$
where 
$$
A_{\sigma}=\{(x_1, x_2,...x_n): x_{\sigma(i)} < x_{\sigma(j)},\, \forall i<j\}.
$$

Now, suppose we have two such points $(x_1, x_2,...x_n)$ and $(x_1,\bar x_2,...,\bar x_n)$ in $S \cap A_{\sigma}$ with \[D_{x_1}c(x_1,x_2,...x_n) =p =D_{x_1}c(x_1,\bar x_2,...,\bar x_n).\]    We claim that either $x_i \geq \bar x_i$ for all $i=2,...n$ or vice versa.  As it is easy to check that $x_i \mapsto D_{x_1}c(x_1,x_2,...,x_n)$  is strictly monotone decreasing on $A_{\sigma}$ for each $i$,  this will establish the result (as if, for example,  $x_i \geq \bar x_i$ for all $i=2,...n$, then the monotonicity implies $D_{x_1}c(x_1,x_2,...x_n)  \leq D_{x_1}c(x_1,\bar x_2,...,\bar x_n)$, with strict inequality as long as
$x_i >\bar x_i$ for at least one $i$).

To see the claim, we use the fact  that the splitting set $S$ is $c$-monotone.  Assume that the claim is false; then there exists non empty disjoint sets $P_+,P_- \subseteq \{1,2,...,n\}$ such that $P_+\cup P_- =\{1,2,...,n\}$, and $x_i \geq \bar x_i$ for $i$ in $P_+$, $x_i \leq \bar x_i$ for $i$ in $P_-$, and at least one of the inequalities is strict in each set.

 For ease of notation, we decompose $(x_1, x_2,...x_n) =(x_+,x_-)$ and $(\bar x_1,\bar x_2,...,\bar x_n) =(\bar x_+,\bar x_-)$, (with the same meaning for $x_+,x_-,\bar x_+$ and $\bar x_-$ as in Definition \ref{def: cmonotonicity}.

By the $c$-montonicity property, we have

$$
c(x_+,x_-) +c(\bar x_+,\bar x_-) \leq c(\bar x_+,x_-) +c( x_+,\bar x_-),
$$
which is equivalent to 

\begin{equation}\label{cmono}
\sum_{i \in P_+, j \in P_-} \frac{1}{|x_i -x_j|}+ \frac{1}{|\bar x_i -\bar x_j|} \leq \sum_{i \in P_+, j \in P_-} \frac{1}{|x_i -\bar x_j|}+ \frac{1}{| \bar x_i - x_j|} .
\end{equation}

Now, for a fixed $i \in P_+, j \in P_-$, consider the paths $x_i(t) = x_i +t(\bar x_i -  x_i )$ and $x_j(s) = x_j +s(\bar x_j -  x_j )$. As both $x$ and $\bar x$ are in the same $A_{\sigma}$ we can assume without loss of generality that $x_i <x_j$ and $\bar x_i <\bar x_j$. Together with the fact that $(\bar x_i-x_j)(\bar x_j-x_j) \leq 0$, (which follows from $i$ being in $P_+$ and $j$ in $P_-$) this implies that for each  point $(t,s) \in [0,1]^2$ we have $x_i(t) < x_j(s)$.  We therefore have 
\begin{eqnarray*}
\frac{\partial ^2}{\partial s\partial t}\Big[\frac{1}{|x_i(t)-x_j(s)|}\Big]&=&\frac{\partial ^2}{\partial s \partial t}\Big[\frac{1}{x_j(s)-x_i(t)}\Big]\\
&=&-2\frac{1}{[x_j(s)-x_i(t)]^3}(\bar x_i -  x_i )(\bar x_j -  x_j )\\
&\geq&0.
\end{eqnarray*}
Now note that this implies 

\begin{eqnarray*}
 \frac{1}{|x_i -x_j|}+ \frac{1}{|\bar x_i -\bar x_j|} -\frac{1}{|x_i -\bar x_j|}- \frac{1}{| x_i -\bar x_j|} &=&\int_0^1\int_0^1 \frac{\partial ^2}{\partial s\partial t}\Big[\frac{1}{|x_i(t)-x_j(s)|}\Big]dsdt\\
&\geq& 0.
\end{eqnarray*}
Furthermore, the inequality is strict for at least one $i \in P_+$ and one $j \in P_-$; summing over $i \in P_+$ and $j \in P_-$,  this violates \eqref{cmono}; this contradiction establishes the claim.

So far, we have proven that  we can have at most one point $(x_2,...x_n) $ such that  $(x_1,x_2,...x_n) \in S \cap A_{\sigma}$ and $D_{x_1}c(x_1,...x_n) =p$ for each permutation $\sigma$.  

Now, note that if $p \geq 0$, we cannot have \textit{any} solutions to $p = D_{x_1}c(x_1,x_2,...x_n) $ in the region where $x_1 <x_i$ for all $i=2,3,...n$ (as clearly  $D_{x_1}c(x_1,x_2,...x_n) <0$ there).  A similar argument applies to the region $x_1 <x_i$ for all $i=2,3,...n$
when $p \leq 0$.  Therefore, we do not have any solutions in $A_{\sigma}$ when $\sigma(1)=1$ in the first case, or when $\sigma(1)=n$ in the second case.  In either case, there are $(n-1)!$ permutations with the appropriate property.   

This means that there are at most $n!-(n-1)!$ permutations $\sigma$ for which we potentially have one solution $(x_2,...x_n) $ to $D_{x_1}c(x_1,...x_n) =p$ on $ S \cap A_{\sigma}$.  We therefore have at most $n!-(n-1)!$ solutions to the equation $p = D_{x_1}c(x_1,x_2,...x_n) $ on the splitting set $S$ .

 \hfill $\square$\\

\begin{corollary} Assume that $c$ is finite $\mu_1\otimes...\otimes \mu_n$ almost everywhere  and there exists  
 a finite
 transport plan $\pi_0$ for $(C)$. If $\mu_1$ is non-atomic and absolutely continuous with respect to the $1-$dimensional Lebesgue measure then there exist  $k \leq n!-(n-1)!,$ a sequence of non-negative measurable real functions 
$\{\alpha_i\}_{i=1}^k $ on $\R$ 
 and,  Borel measurable maps $G_1,...,G_k: \R \to  \R^{n-1}$ such that
\begin{equation}\label{late0}\pi_0 =\sum_{i=1}^k \alpha_i (\text{Id} \times  G_i)_\# \mu_1, \qquad \end{equation}
where  $\sum_{i=1}^k \alpha_i(x) =1$ for $\mu_1$-a.e. $x \in \R$.
\end{corollary}
\textbf{Proof.} By a similar argument as in the proof of Theorem (4) in \cite{B-D-G}, there exist potentials $u_1,..., u_n$ with $u_i\in \mathcal{L}^1(\mu_i)$, taking values in $[-\infty, \infty)$ and $ c(x_1,..., x_n) \geq \sum_{i=1}^nu(x_i)$
such that \[\int_{\R^n} c(x_1,...x_n)\, d\pi_0=\sum_{i=1}^n \int_\R u(x_i) \, d\mu_i,\]
and $u_1$ is differentiable $\mu_1$ almost everywhere. It also follows from Proposition \ref{col} that 
$c$ is $n!-(n-1)!$ twisted in splitting sets. Thus, the result follows from Theorem \ref{main}.
\hfill $\square$\\

As shown in \cite{cdpdm}, when the marginals are all equal (the physically relevant case) the solution is concentrated on $(n-1)!$ graphs; our result guarantees it is concentrated on \textit{at most} $n!-(n-1)!$, even when the marginals differ.  It is unclear to us whether this is sharp; that is, whether there really are solutions which concentrate on $n!-(n-1)!$ graphs.  If so, this can only happen when the marginals are not all equal.\\

\subsection{A cost satisfying local $1-$twistedness}
We exhibit now an example of a multi-dimensional cost function  which is locally (but not globally) $1-$twisted on splitting sets, together with a solution to the optimal transport problem which concentrates on several (rather than one) graphs.

\begin{proposition}
Take, for $x=(x^1,x^2),y=(y^1,y^2),z=(z^1,z^2) \in \mathbb{R}^2$,
\begin{equation}\label{examplecost}
c(x,y,z)=-e^{x^1+y^1}\cos(x^2-y^2)-e^{x^1+z^1}\cos(x^2-z^2)-e^{y^1+z^1}\cos(z^2-y^2).
\end{equation}
Then $c$ is locally twisted on splitting sets.
\end{proposition} 
\textbf{Proof.}
We first recall that, for three marginal costs, the condition
$$
D^2_{xy}c[D^2_{zy}c]^{-1}D^2_{zx}c <0
$$
is equivalent to the matrix $g$ in \eqref{gmatrix} having signature $(2k,2,0)$ (see \cite{P4});  hence, by Proposition \ref{gsignature}, this condition implies the local twist on splitting sets property).  We will  show that this condition holds for the cost \eqref{examplecost}.
Note that 

\begin{equation*} \qquad
D^2_{xy}c=
-\begin{bmatrix}
 e^{x^1+y^1}\cos(x^2-y^2) & e^{x^1+y^1}\sin(x^2-y^2) \\
-e^{x^1+y^1}\sin(x^2-y^2) & e^{x^1+y^1}\cos(x^2-y^2) \\
\end{bmatrix}.
\end{equation*}
Up to multiplication by $-e^{x^1+y^1}$, this is a rotation matrix through the angle $x^2-y^2$.  Similarly, $ D^2_{yz}c$ is $-e^{y^1+z^1}$ multiplied by a rotation through  $y^2-z^2$, so its inverse transpose, $[D^2_{zy}c]^{-1}$, is $-\frac{1}{e^{y^1+z^1}}$ multiplied by the rotation through $y^2-z^2$ (recall that, for a rotation matrix $A$, $A^T=A^{-1}$).  Finally, $ D^2_{zx}c$ is $-e^{x^1+z^1}$ multiplied by a rotation through  $z^2-x^2$.  Therefore, the product $D^2_{xy}c[D^2_{zy}c]^{-1}D^2_{zx}c$ is, up to a negative multiplicative constant, rotation through $x^2-y^2 +y^2-z^2 +z^2-x^2=0$; that is, the product is a multiple of the identity:

$$
D^2_{xy}c[D^2_{zy}c]^{-1}D^2c_{zx}c=-e^{x^1+y^1}\frac{1}{e^{y^1+z^1}}e^{x^1+z^1} I <0.
$$

Therefore, by Proposition \ref{gsignature}, this cost is locally twisted on splitting sets.\hfill $\square$\\

In fact, for this cost, we can exhibit explicitly an optimizer which  concentrates on several graphs.  To this end,  note that 

$$
-e^{x^1+y^1}\cos(x^2-y^2) \geq -e^{x^1+y^1} \geq -e^{2x^1}/2 -e^{2y^1}/2
$$
and we have equality only when $x^1=y^1$ and $x^2-y^2$ is an integer multiple of $2\pi$. Applying similar reasoning to the other pieces of $c$, we get
$$
c(x,y,z) \geq -e^{2x^1} -e^{2y^1}-e^{2z^1}
$$
with equality only when 
$$
(x,y,z) \in S=\{(x,y,z): x^1=y^1=z^1\text{ and } x^2-y^2 = 2k\pi, \text{ } x^2-z^2 =2l\pi\text{ for some integers }k,l \}
$$
This easily implies that any measure $\gamma$ concentrated on the set $S$ is optimal in $(MK)$ for its marginals,  as  $-e^{2x^1}, -e^{2y^1}$ and $-e^{2z^1}$ serve as Kantorovich potentials.

Locally, the set $S$ looks like a graph, but globally, it is the union of countably many graphs (or finitely many, if we restrict to compact domains).

\section{Uniqueness issues}

In this section, we consider the  uniqueness of solutions concentrated on several graphs. We begin with a criterion for extremality of measures of type (\ref{late}) in $\Pi(\mu_1,...,\mu_n)$.  Note that, as the Kantorovich problem is a minimization of a linear functional over a convex set $\Pi(\mu_1,...,\mu_n)$, it necessarily has at least one solution which is extremal in that set.  Extremality of a solution $\gamma \in \Pi(\mu_1,...,\mu_n)$ to $(MK) $ is therefore a {\em necessary} condition for $\gamma$ to be the unique  solution.

\begin{theorem}\label{uniq} Let $X_1, X_2, X_3$  be  Polish spaces equipped with  Borel probability measures $\mu_i$ on $X_i$,
 and    let  $\{G_i=(H_i,K_i)\}_{i=1}^k$ be a finite  sequence of measurable maps from $ X_1$ to $X_2 \times X_3$.   
Assume that   the following assertions hold:
\begin{enumerate}
\item[(i)] For each $i \geq 1 $ the map  $H_i$ is injective and onto.
\item [(ii)] For each $j \geq 2$ the map $K_j$ is injective and $Ran(K_i)\cap Ran(K_j)=\emptyset$ for all $ i\geq 2$ with $i \not=j.$ 
\item [(iii)] There exists a bounded measurable function $\theta:X_3 \to \R$ with the property that for each $i \geq 2,$  \[\theta\big (K_1 \circ H^{-1}_1(x_2)\big )>\theta\big (K_i \circ H^{-1}_i(x_2)\big ) \qquad \forall x_2 \in Dom\big (K_1 \circ H^{-1}_1 \big )\cap Dom\big (K_i \circ H^{-1}_i \big ).\]
\end{enumerate}
Then each  $\gamma \in \Pi(\mu_1, \mu_2, \mu_3)$ that is concentrated on $ \cup_{i=1}^k Graph (G_i)$ is an extremal point of $\Pi(\mu_1, \mu_2, \mu_3).$
 \end{theorem}

The proof of the preceding theorem is fairly long and so is deferred to the end of this section (see subsection \ref{sect: proof uniq}).  We note that we have stated the result for $n=3$ only to keep the presentation as simple as possible;  it has a straight forward generalization to larger $n$.  

The theorem may be used in certain circumstances to deduce uniqueness of the optimizer; in particular, if one can show that \emph{all} solutions  must be concentrated on a single collection of several graphs, satisfying the hypotheses in the Theorem, then uniqueness follows immediately.  The following example illustrates this concept.

To the best of our knowledge, this is the first example of a  solution to a multi-marginal optimal transport problem which is both  1) unique and  2) not concentrated on a single graph.  Although the example is admittedly somewhat artificial, we still believe it is of interest as this type of behavior has not been observed before in multi-marginal problems.
\begin{example}\textbf{A unique solution concentrated on two graphs} 

Let $c(x,y,z) =(x-y)^2+(x-z)^2(x-z+1/2)^2$ on $X_1 \times X_2 \times X_3=[0,1]\times[0,1] \times [0,3/2]$. Consider the maps $(H_1,K_1):X_1 \rightarrow X_2 \times X_3$ given by $(H_1,K_1)(x) =(x,x)$ and   $(H_2,K_2):X_1 \rightarrow X_2 \times X_3$ given by $(H_2,K_2)(x) =(x,x+\frac{1}{2})$.  We note that

$$
c(x,y,z) \geq 0
$$
with equality if and only if $(y,z) =(H_i,K_i)(x)$, for $i=1$ or $2$.  Therefore, letting  $\gamma$ be uniform measure (renormalized to have total mass $1$) on the union $S=graph(H_1,K_1) \cup graph(H_2,K_2)$ of these two graphs, we have that $\gamma$ is optimal for it's marginals $\mu_1,\mu_2,\mu_3$.  To see that $\gamma$ is in fact the unique optimizer, assume that $\bar \gamma$ is any other optimizer; it must necessarily concentrate on the same set $S$.  The interpolant $\frac{1}{2}(\gamma +\bar \gamma)$ then also concentrates on $S$ and must the  be optimal too.  However, it is easy to see that $(H_i,K_i)$ satisfy the conditions in Theorem \ref{uniq} (taking, for instance, $\theta(z) =-z$), and so $\frac{1}{2}(\gamma +\bar \gamma)$ must be extremal in $\Pi(\mu_1,\mu_2,\mu_3)$, a contradiction as it is the average of two distinct measures $\gamma, \bar \gamma \in \Pi(\mu_1,\mu_2,\mu_3)$. 
Therefore, $\gamma$ is the unique optimal measure, as desired.\\
\end{example}
We now record a result that shows that, under the $m$-twisted on splitting sets criterion, uniqueness can only fail by reorganizing the mass within a single collection of $m$ graphs.

\begin{theorem} \label{us} Suppose that  the continuously differentiable function  $c$   satisfies  the $m$-twist condition on splitting sets, 
$\mu_1$ is non-atomic  and  any  function $u_1(x_1)$ of the form \eqref{eqn: c-conjugate} which is not identically infinite is differentiable $\mu_1$-almost surely on its domain. Assume that  $\{G_i\}_{i=1}^m$ is  a sequence of measurable functions from $X_1$ to $X_2 \times ...\times X_n$ such that  for each $i\not=j$ the set $\{x; G_i(x)=G_j(x)\}$ is $\mu_1$-negligible.
Let
$\bar \gamma$ be an optimal plan for $(MK)$ such that
\begin{eqnarray}
\bar \gamma =\sum_{k=1}^m \alpha_k (\text{Id}\times  G_k)_\# \mu, \qquad  \big (\alpha_i(x_1) \geq 0  \text{ and } \, \alpha_1(x_1) \alpha_2(x_1)...\alpha_m(x_1)\not=0 \text{ for }  \mu_1-\text{a.e. } \, x_1 \in X_1\big ).
\end{eqnarray}
  Then for any other optimal plan $\gamma$  we have
 \[\text{Supp} (\gamma) \subseteq \text{Supp}(\bar \gamma).\]
\end{theorem}

\textbf{Proof.}
By Kantorovich duality (see \cite{Ke, C-N} for a proof  in the multi-marginal case), there exist   functions $\phi_i \in \mathcal{L}^1(\mu_i),$ $i=1,...,n$  such that 
\[\phi_i(x_i)=\inf_{x_j \in X_j, \, j \not=i} \Big\{c(x_1,...,x_n)-\sum_{j\not=i} \phi_j(x_j) \Big\},\]
for all $ x_i \in X_i$ and
 \begin{equation*}
 \int c \, d\bar \gamma= \sum_{i=1}^n\int_{X_i} \phi_i(x_i) \, d \mu_i.
 \end{equation*}
Let $S$ be the $c$-splitting set generated by the $n$-tuple $(\phi_1,...,\phi_n)$, that is, \[S=\Big \{(x_1,...,x_n); \, c(x_1,..., x_n)=\sum_{i=1}^n\phi_i(x_i)  \Big \}.
\]
 It follows that
\[ \int c(x_1,...,x_n) \, d \bar  \gamma= \int  \sum_{i=1}^n \phi_i(x_i)\, d\bar \gamma,  \]
from which we obtain
\[\sum_{k=1}^m \int_{X_1} \alpha_i(x_1) c(x_1,G_i x_1) \, d \mu_1=\sum_{k=1}^m \int_{X_1} \alpha_i(x_1) \Big[ \phi_1(x_1)+\psi (G_i x_1)\Big] \, d \mu_1,\]
where $\psi(x_2,...,x_n)= \sum_{i=2}^n \phi_i(x_i).$
It then follows that
\[\sum_{k=1}^m \int_{X_1} \alpha_i(x_1) \Big [c(x_1,G_i x_1)-\phi_1(x_1)-\psi (G_i x_1)\Big] \, d \mu_1=0.\]
Since each integrand in the latter expression is non-negative  and $\mu_1$ almost surely $\alpha_1(x_1) \alpha_2(x_1)...\alpha_m(x_1)\not=0,$  we have  that
\[c(x_1,G_i x_1)=\phi(x_1)+\psi (G_i x_1)\qquad \mu_1-a.e. \quad \forall i \in\{1,...,m\}.\]
Consequently we  obtain,
\begin{equation}\label{uniq0}D_{x_1} c(x_1,G_i x_1)=D \phi_1(x_1)\qquad \mu_1-a.e. \quad \forall i \in\{1,...,m\}.\end{equation}
Note also that  for  $i\not=j$ the set $\{x_1 \in X_1; \, \, G_i(x_1)=G_j(x_1)\}$ is a null set with respect to the measure $\mu_1.$ This together with (\ref{uniq0}) and
 the $m$-twist
condition on $S$ imply that  the cardinality of
 the set $\{G_1x_1,...,G_m x_1\}$ is $m$ for $\mu_1$-a.e. $x_1 \in X_1.$\\
 Now assume that  $\gamma$ is also  an optimal plan
   of  $(MK).$   It follows from Theorem \ref{main} that
 there exist  a sequence of non-negative functions
$\{\beta_i\}_{i=1}^m $
 and,  Borel measurable maps $T_1,...,T_m: X_1 \to X_2 \times...\times X_n$ such that
\begin{eqnarray*}
\gamma =\sum_{i=1}^m \beta_i (\text{Id} \times  T_i)_\# \mu_1.
\end{eqnarray*}
By a similar argument as above one obtains
\begin{equation*}\beta_i(x)\big [D_{x_1} c(x_1,T_i x)-D\phi_1(x_1)\big]=0\qquad \mu_1-a.e. \quad \forall i \in\{1,...,m\}.\end{equation*}
For each $i$ define $\Omega_i=\{x_1 \in X_1; \beta_i(x_1)\not=0\}.$
Since the cardinality of the set $\{G_1x_1,...,G_m x_1\}$ is $m$ for $\mu_1$-a.e. $x_1 \in X_1$ and since  $c$ satisfies the $m$-twist condition on $S$
we have  that for each $i,$ $\{T_i x_1\} \subseteq \{G_1x_1,...,G_m x_1\}$ for $\mu_1$-a.e. $x_1 \in \Omega_i$. This completes the proof. \hfill $\square$\\

\subsection{Permutation symmetric problems}
Finally, we turn our attention to permutation symmetric problems; that is, problems for which the cost is symmetric under permutation of its arguments, and the marginal are all the same.  As discussed in the introduction, these problems have important applications in physics and economics. 

We first demonstrate with an example that it is in fact possible to have unique solutions to these problems concentrated on several graphs.  In \cite{P3}, it was noted that a problem of this type could not admit a unique solution concentrated on a single graph, $\{(x,F_2(x),...,F_n(x)  \}$, unless each component of that  graph was the identity almost everywhere, $F_i(x) =x$ for a.e. $x$, for each $i=2,...,m$.  The example below illustrates that one can have somewhat less trivial unique solutions if we relax the structural requirements to being concentrated on two graphs instead of one.  However, as Proposition \ref{symm_uniqueness} below implies, unique optimizers are still fairly special;  this proposition asserts that one can not have unique solutions, except under very particular conditions.  

 \begin{example}\textbf{(Another unique solution concentrated on two graphs)}\\
Take $c(x,y,z) =xyz$ on $[-1,1]^3$, with each marginal equal to

$$
\mu=\frac{1}{3}[\mathcal{L}_{[-1,0]} +2\mathcal{L}_{[0,1]}]
$$
Define $G_1:[-1,1] \rightarrow [-1,1]^2$ by $G_1(x)=(-x,|x|)$ and $G_2:[-1,1] \rightarrow [-1,1]^2$ by $G_2(x) =(x,-|x|)$.  
Note that by the geometric-arithmetic mean inequality, we have

$$
|xyz| \leq \frac{|x|^3}{3}+ \frac{|y|^3}{3}+ \frac{|z|^3}{3},
$$
with equality only when $|x| =|y| =|z|$ and so
$$
c(x,y,z) =xyz \geq -|xyz| \geq -\frac{|x|^3}{3}-\frac{|y|^3}{3}-\frac{|z|^3}{3}
$$
with equality only when either one or three of $x,y$ and $z$ are non-positive, and $|x| =|y|=|z|$; that is, we have equality precisely on the graphs of $G_1$ and $G_2$.  Therefore, if we can find a $\gamma \in \Pi(\mu,\mu,\mu)$ concentrated on these two graphs, it is optimal; if there is a unique such $\gamma$, then the optimal measure is unique (as any probability measure concentrated outside these sets has larger total cost).

Now, finding a measure $\gamma \in \Pi(\mu,\mu,\mu)$ of the form

\begin{equation}\label{eqn: optimal measure}
\gamma = \alpha_1(x) (Id \times G_1)_\#\mu +  \alpha_2(x)(Id \times G_2)_\#\mu 
\end{equation}
amounts  to solving, for almost every $x\geq 0$, the system of equations:

\begin{eqnarray*}
\alpha_1(x) +\alpha_2(x) &=& 1\\
\alpha_1(-x) +\alpha_2(-x) &=& 1\\
\frac{1}{3}\alpha_1(-x) +\frac{2}{3}\alpha_2(x) &=& \frac{2}{3}\\
\frac{2}{3}\alpha_1(x) +\frac{1}{3}\alpha_2(-x) &=& \frac{1}{3}\\
\frac{2}{3}\alpha_1(x) +\frac{1}{3}\alpha_1(-x) &=& \frac{2}{3}\\
\frac{2}{3}\alpha_2(x) +\frac{1}{3}\alpha_2(-x) &=& \frac{1}{3}.
\end{eqnarray*}
A straightforward calculation shows that the unique solution to these equations is given by
$$
\alpha_1(x)= \left\{
  \begin{array}{l l}
    \frac{1}{2} & x>0\\
    1&x<0,
  \end{array} \right.
$$
and
$$
\alpha_2(x)= \left\{
  \begin{array}{l l}
    \frac{1}{2} & x>0\\
    0 &x<0.
  \end{array} \right.
$$
The measure given by equation  \eqref{eqn: optimal measure}, with these coefficients, is therefore the unique solution.
\end{example}

The preceding example demonstrates that it is possible to have a unique solution to a symmetric problem which concentrates on two graphs.  The components of these graphs need not be the identity, but note that for each $G_i =(H_i,K_i)$ and each $x$, one of the following holds:

\begin{enumerate}
\item $H_i(x)=x$,
\item $K_i(x)=x$, or, 
\item $H_i(x)=K_i(x)$. 
\end{enumerate}
 
 In particular, $\gamma(S_1 \times S_2 \times S_3) =0$ whenever the sets $S_1,S_2,S_3 \subseteq [-1,1]$ are pairwise disjoint.  This is essentially the only way that multiple maps can support a unique solution to a symmetric problem, as the following result shows.

\begin{proposition}\label{symm_uniqueness}
Let $n=3$.  Suppose that the spaces $X_i :=X$ and marginals $\mu_i:=\mu$ are the same for each $i$, and that $c(x,y,z)$ is symmetric with respect to any permutation of the arguments.  Assume that there exist mutually disjoint sets $S_1,S_2,S_3 \subseteq X$, and an optimal measure $\gamma$ that charges $S:=S_1 \times S_2 \times S_3$; that is, $\gamma(S) >0$.  Then the solution is nonunique.
\end{proposition}
As should be clear from the proof, a similar result holds for $n \geq 4$.

\textbf{Proof.} The proof is by contradiction; assume that $\gamma$ is the unique solution.  Letting $\sigma$ be any permutation on the variables $(x,y,z)$, this uniqueness implies that $\sigma_{\#}\gamma =\gamma$.
  
Now, note that

$$
\sigma _\#(\gamma_{S_1 \times S_2 \times S_3}) =\gamma_{\sigma(S_1 \times S_2 \times S_3)},
$$ 
where $\gamma_{S_1 \times S_2 \times S_3}$ denotes the restriction of $\gamma$ to the set $S_1 \times S_2 \times S_3$.
Therefore, by the symmetry of $c$, we have
\begin{equation}\label{measure_symm}
\int_{X^3} c(x,y,z)d\gamma_{S_1 \times S_2 \times S_3}= \int_{X^3}c(x,y,z)d(\sigma _\#(\gamma_{S_1 \times S_2 \times S_3})) =\int_{X^3}c(x,y,z)d\gamma_{\sigma(S_1 \times S_2 \times S_3)}.
\end{equation}
It is clear that each $\gamma_{\sigma(S_1 \times S_2 \times S_3)} \leq \gamma$, and that the measures 
$$
\frac{1}{6}\sum_{\sigma \in \mathcal{S}_3}\gamma_{\sigma(S_1 \times S_2 \times S_3)}
$$
and
$$
\frac{1}{3}[\gamma_{S_1 \times S_2 \times S_3} +\gamma_{S_3 \times S_1 \times S_2}+\gamma_{S_2 \times S_3 \times S_1}]
$$
share the same marginals, where  $\mathcal{S}_3$ is the set of permutations on $3$ letters.  It then follows that
$$
\bar \gamma =\gamma +\frac{1}{6}\sum_{\sigma \in S_3}\gamma_{\sigma(S_1 \times S_2 \times S_3)}-\frac{1}{3}[\gamma_{S_1 \times S_2 \times S_3} +\gamma_{S_3 \times S_1 \times S_2}+\gamma_{S_2 \times S_3 \times S_1}]
$$
is a nonnegative measure with the same marginals as $\gamma$.  On the other hand, $\bar \gamma \neq \gamma$, as 
$$
\frac{1}{6}\sum_{\sigma \in S_3}\gamma_{\sigma(S_1 \times S_2 \times S_3)} \neq \frac{1}{3}[\gamma_{S_1 \times S_2 \times S_3} +\gamma_{S_3 \times S_1 \times S_2}+\gamma_{S_2 \times S_3 \times S_1}].
$$
Finally,  as by \eqref{measure_symm},

$$
\int_{X^3} c(x,y,z) d\gamma =\int_{X^3} c(x,y,z)d\bar \gamma,
$$
$\bar \gamma$ is also optimal, which contradicts uniqueness. \hfill $\square$\\

\subsection{Extremal  measures with fixed marginals}\label{sect: proof uniq}
We shall now prove  Theorem \ref{uniq}; the proof requires a few preliminary results.
 We begin with  the following  result  from \cite{Moa3}.

\begin{theorem}\label{aper} Let $X$ and  $Y$  be  Polish spaces equipped with  Borel probability measures $\mu$ on $X$
 and $\nu$ on $Y,$ and    let  $\{T_i\}_{i=1}^k$ be a (possibly infinite) sequence of measurable maps from $ X$ to $Y$.    
Assume that   the following assertions hold:
\begin{enumerate}
\item [(i)]For each $i \geq 2 $ the map  $T_i$ is injective on the set 
\[D_i:=\big \{x \in Dom(T_1)\cap Dom(T_i); \,\, \, T_1x\not=  T_i x\big \},\] and 
$Ran(T_i)\cap Ran(T_j)=\emptyset$
for all $ i,j\geq 2$ with $i \not=j.$
\item [(ii)] There exists a bounded measurable function $\theta:Y \to \R$ with the property that $\theta(T_1x)>\theta(T_ix)$ on $D_i.$
\end{enumerate}
Then there exists at most one $\gamma \in \Pi(\mu, \nu)$ that is concentrated on $ \cup_{i=1}^k Graph (T_i).$ Moreover,  any $\gamma \in \Pi(\mu, \nu)$ that is concentrated on $ \cup_{i=1}^k Graph (T_i)$ is an extremal point of $\Pi(\mu, \nu).$
 \end{theorem}

Our proof exploits this result, together with a connection between extremality in the two and three marginal cases, established in the following lemma.
\begin{lemma}\label{trip} Let $\gamma$ be a  measure in $\Pi(\mu_1, \mu_2, \mu_3)$ and let $\nu$ be the projection of $\gamma$ onto $X_2 \times X_3.$  If $\gamma$ is an extremal point of $\Pi(\mu_1, \nu)$ and $\nu$ is an extremal point of $\Pi(\mu_2, \mu_3)$ then $\gamma$ is an extremal point of $\Pi(\mu_1, \mu_2, \mu_3)$.
\end{lemma}

\textbf{Proof.}  Let $\gamma_1, \gamma_2 \in \Pi(\mu_1, \mu_2, \mu_3)$ and $0<t<1$ be such that $\gamma=t\gamma_1+(1-t)\gamma_2.$ For $i=1,2$ denote by $\nu_i$ the projection of $\gamma_i$ on $X_2 \times X_3.$ It is obvious that $\nu_1, \nu_2 \in \Pi(\mu_2, \mu_3),$  and that $\nu=t\nu_1+(1-t)\nu_2.$ Since $\nu $ is an extremal point of $\Pi(\mu_2, \mu_3)$ and $0<t<1$ we must have $\nu=\nu_1=\nu_2,$ which meands that $\gamma_1,\gamma_2 \in \Pi(\mu_1, \nu)$.  On the other hand $\gamma$ is an extremal point of $\Pi(\mu_1,\nu)$; together with $\gamma=t\gamma_1+(1-t)\gamma_2$, this implies $\gamma=\gamma_1=\gamma_2,$ from which we obtain the desired conclusion. \hfill $\square$\\

Let $(X, \Sigma_1,\mu) $ and $(Y, \Sigma_2, \nu)$  be  two measure spaces.   A map   $T: X \to Y$ is said to be a  push forward from $\mu$ to $\nu$, dented by
 $T_\# \mu=\nu,$ if 
\[\forall f \in {\mathcal L}^1(\nu), \quad \, f \circ T \in {\mathcal L}^1(\mu) \, \& \, \int_{X}f(Tx) \, d\mu=\int_{Y}f(y) \, d\nu,\]
where  ${\mathcal L}^1(\nu)$ is the set of $\nu-$integrable functions on $Y.$\\

\textbf{Proof of Theorem \ref{uniq}.} 
Let $\gamma \in  \Pi(\mu_1, \mu_2, \mu_3)$ be a triply stochastic measure that is concentrated on $\cup_{i=1}^k Graph(G_i).$ By definition there exist non-negative Borel measurable real functions  $\alpha_i: X_1 \to \R$  with $\sum_{i=1}^k \alpha_i(x_1)=1$  such that $\gamma=\sum_{i=1}^k \alpha_i(Id \times G_i)_\# \mu_1,$ i.e.,
\[\int f\, d\gamma=\sum_{i=1}^k \int_{X_1} \alpha_i(x_1) f\big(x_1, H_i x_1, K_i x_1 \big) \, d\mu_1.\]

  Denote by $\nu$ the projection of $\gamma$ on $X_2 \times X_3.$ 
We proceed with the proof in several steps.\\

{\it Step 1.} In this step we show that the projection $\nu \in \Pi(\mu_2, \mu_3)$ of $\gamma$ onto $X_2 \times X_3$  is concentrated on $\cup_{i=1}^k Graph (K_i \circ H_i^{-1}).$  For $1 \leq i \leq k,$ define the measures $\mu_{i,1}$ and $\mu_{2,i}$ by $d \mu_{i,1}=\alpha_i d\mu_1$ and  $\mu_{2,i}={H_i}_\# \mu_{1,i}.$  For every measurable set $B \in \mathcal{B}(X_2)$ we have 
\begin{eqnarray*}
\mu_2(B)=\gamma (X_1 \times B \times X_3)=\sum_{i=1}^k\int_{X_1} \alpha_i(x_1)\chi_{B}  ( H_i x_1) \, d\mu_1&=& \sum_{i=1}^k\int_{X_1} \chi_{B}  ( H_i x_1) \, d\mu_{1,i}\\
&=&\sum_{i=1}^k\int_{X_2}\chi_{B}  (x_2) \, d\mu_{2,i}\\
&=&\sum_{i=1}^k\int_B \, d\mu_{2,i}=\sum_{i=1}^k\mu_{2,i}(B)
\end{eqnarray*}
Thus, for all $B \in \mathcal{B}(X_2)$ we have $\mu_2(B)=\sum_{i=1}^k\mu_{2,i}(B)$. It implies that each $\mu_{2,i}$ is absolutely continuous with respect to $\mu_2.$ Therefore, for each $i$ there exists a nonnegative real function $\beta_i: X_2 \to \R$ such that $d\mu_{2,i}=\beta_i d \mu_2.$  This together with $\mu_2(B)=\sum_{i=1}^k\mu_{2,i}(B)$ yield that 
\[\mu_2(B)=\sum_{i=1}^k\int_B  \beta_i(x_2) \, d\mu_2, \qquad \forall B \in \mathcal{B}(X_2),\]
from which we obtain $\sum_{i=1}^k  \beta_i(x_2)=1$ for $\mu_2$ almost every $x_2 \in X_2.$
Take a measurable bounded function $f: X_2 \times X_3 \to \R.$
It follows that
\begin{eqnarray*}
\int_{ X_2 \times X_3} f(x_2,x_3)\, d\nu=\int_{X_1\times X_2 \times X_3} f(x_2,x_3)\, d\gamma= &=&\sum_{i=1}^k\int_{X_1} \alpha_i(x_1)f ( H_i x_1, K_i x_1) \, d\mu_1\\
&=&\sum_{i=1}^k\int_{X_1} f ( H_i x_1, K_i x_1) \, d\mu_{1,i}\\
&=&\sum_{i=1}^k\int_{X_1} f (x_2, K_i\circ H_i^{-1}x_2) \, d\mu_{2,i}\\
&=&\sum_{i=1}^k\int_{X_2}\beta_i(x_2)f (x_2, K_i\circ H_i^{-1}x_2) \, d\mu_2.
\end{eqnarray*}
This implies that $\nu=\sum_{i=1}^k\beta_i (Id \times K_i\circ H_i^{-1})_\#{\mu_2}.$\\

{\it Step 2.} In this step we show that $\nu$ is an extremal point of $\Pi(\mu_2, \mu_3).$ We shall make use of Theorem \ref{aper} to prove this. By assumption we have that $K_i$ is injective and $Ran(K_i)\cap Ran(K_j)=\emptyset$ for all $i,j\geq 2$ with $i\not=j.$  This implies  that $K_i\circ H_i^{-1}$ is injective for all $i\geq 2$ and $Ran(K_i\circ H_i^{-1})\cap Ran(K_j\circ H_j^{-1})=\emptyset$
from which assumption $(i)$ of Theorem \ref{aper} follows. The second assumption in Theorem \ref{aper} also follows due to hypothesis $(iii)$ in Theorem \ref{uniq} and therefore the extremality of $\nu$ in $\Pi(\mu_2, \mu_3)$ follows.\\

{\it Step 3.} We show that $\gamma$ is an extremal point of $\Pi(\mu_1, \nu).$ We shall again use Theorem \ref{aper} to prove this part. Note that $G_i=(H_i, K_i)$ is injective for $i\geq 2$
as both $H_i$ and $K_i$ are injective. For $i,j \geq 2$  with $i\not=j$ we have $Ran(G_i)\cap Ran(G_j)=\emptyset $ as by assumption $Ran(K_i)\cap Ran(K_j)=\emptyset.$   Define $\phi: X_2 \times X_3 \to \R$
by
\begin{eqnarray*}
		\phi(x_2,x_3)=\left\{
		\begin{array}{ll}
			(k-1)\theta(x_3)-\sum_{i=2}^k \theta (K_i\circ H_i^{-1} x_2), & x_3 \in X_3\setminus \cup_{i=2}^k Ran(K_i),\\
			
			\theta(x_3)-\theta (K_j\circ H_j^{-1}x_2), &   x_3 \in Ran(K_j) \text{ for some } j\geq 2.\\
		\end{array}
		\right.
	\end{eqnarray*}
Since $K_i$ is injective for $i\geq 2$ it follows that $Ran(K_i)$ is Borel measurable (\cite{Bo}, Theorem 6.8.6). Thus, $\phi$ is a bounded Borel measurable function.
We show that $\phi(G_1 (x_1)) >\phi(G_j(x_1))$ for all $j\geq 2$ and $x_1 \in Dom(G_1) \cap Dom(G_j).$  Fix 
$x_1 \in Dom(G_1) \cap Dom(G_j).$ If $K_1 x_1 \in X_3\setminus \cup_{i=2}^k Rang(K_i)$ it follows that 
\begin{eqnarray*}
\phi(G_1 (x_1))=\phi(H_1 x_1, K_1 x_1)&=&(k-1)\theta (K_1 x_1)-\sum_{i=2}^k \theta (K_i\circ H_i^{-1}\circ H_1 x_1)\\&>&
  (k-1)\theta (K_1 x_1)-\sum_{i=2}^k \theta (K_1\circ H_1^{-1}\circ H_1 x_1) \\
  &=&(k-1)\theta (K_1 x_1)-(k-1)\theta (K_1 x_1)=0,
\end{eqnarray*}
where in the second line we have used the fact that $\theta \circ K_i\circ H_i^{-1} < \theta \circ K_1 \circ H_1^{-1}.$ 
If now  $K_1 x_1 \in Ran(K_j)$ for some $j \geq 2$ we obtain
\begin{eqnarray*}\phi(G_1 (x_1))=\phi(H_1 x_1, K_1 x_1)&=&\theta(K_1 x_1)-\theta (K_j\circ H_j^{-1}\circ H_1 x_1)\\
&>& \theta(K_1 x_1)-\theta (K_1\circ H_1^{-1}\circ H_1 x_1)\\&=& \theta(K_1 x_1)-\theta(K_1 x_1)=0.\end{eqnarray*}
It then follows that in both cases we have $\phi(G_1 (x_1))>0.$
We now sow that  $\phi(G_j(x_1))=0$. In fact,
\[\phi(G_j(x_1))=\phi(H_j x_1, K_j x_1)=\theta(K_j x_1)-\theta(K_j\circ H_j^{-1}\circ H_j x_1)=0.\] 
Therefore for each $x_1 \in Dom(G_1) \cap Dom(G_j)$ we have
\[\phi(G_1 (x_1))>0=\phi(G_j(x_1)).\]
It now follows from Theorem \ref{aper} that $\gamma$ is an extremal point of $\Pi(\mu_1, \nu).$
Therefore, it follows from Lemma \ref{trip} together with steps 2 and 3 that $\gamma$ is an extremal point of $\Pi(\mu_1, \mu_2, \mu_3).$ \hfill $\square$\\

\end{document}